\documentclass{amsart}

\usepackage{etex}
\usepackage{amsmath, amssymb}
\usepackage{mathtools}
\usepackage[frame,cmtip,arrow,matrix,line,graph,curve]{xy}
\usepackage{graphpap, color, paralist, pstricks}
\usepackage[mathscr]{eucal}
\usepackage{mathabx}
\usepackage[pdftex,colorlinks,backref=page,citecolor=blue]{hyperref}
\usepackage{tikz}

\usepackage{amsmath}
\usepackage{amsthm,amsfonts,amssymb,mathrsfs}
\usepackage{epic,eepic}
\usepackage{yfonts}
\usepackage{paralist,enumerate}

\theoremstyle{definition}
\newtheorem{theorem}{Theorem}
\newtheorem{definition}[theorem]{Definition}

\theoremstyle{remark}

\newtheorem{example}[theorem]{Example}

\begin{document}

\title{Combinatorial Applications of the Hodge-Riemann Relations}

\author{June Huh}

\address{Institute for Advanced Study, Fuld Hall, 1 Einstein Drive, Princeton, NJ, USA.}
\email{junehuh@ias.edu}

\maketitle

\section{Log-concave and unimodal sequences}

Logarithmic concavity is a property of a sequence of  real numbers, occurring throughout  algebra, geometry, and combinatorics.  A sequence of real numbers $a_0,\ldots,a_{d}$ is \emph{log-concave} if
\[
a_i^2\ \ge\ a_{i-1}a_{i+1} \ \ \text{for all $i$}.
\]
When all the entries are positive, the log-concavity implies unimodality, a property easier to visualize: the sequence  is  {\em unimodal} if there is an index $i$ such that
\[
a_0\leq\cdots \leq a_{i-1} \leq a_i \geq a_{i+1} \geq\cdots\geq a_{d}.
\]
A rich variety of log-concave and unimodal sequences arising in combinatorics can be found in the surveys \cite{Brenti94,Stanley89,Stanley00}.
For an extensive discussion of log-concavity and its applications in probability and statistics, see \cite{DJ88,MOA11,SW14}.

Why do  natural and interesting sequences often turn out to be log-concave?
Below we give one of many possible explanations, from the viewpoint of ``standard conjectures''.
To illustrate, we  discuss three combinatorial sequences appearing in \cite[Problem 25]{Stanley00},  in Sections \ref{Chromatic}, \ref{Independent},  and \ref{Open}.
Another heuristic, based on the physical principle  that the entropy of a system should be concave as a function of the energy, can be found in \cite{Okounkov03}.

Let $X$ be a mathematical object of ``dimension'' $d$.
Often it is possible to construct from $X$ in a natural way
 a graded vector space over the real numbers
\[
A(X)=\bigoplus_{q=0}^d A^q(X),
\]
a symmetric bilinear map 
$
\mathrm{P}:A(X) \times A(X) \rightarrow \mathbb{R},
$
and a graded linear map 
$\mathrm{L}:A^\bullet(X) \rightarrow A^{\bullet+1}(X)$ that is symmetric with respect to $\mathrm{P}$.
The linear operator $\mathrm{L}$ usually comes in as  a member of a family $\mathrm{K}(X)$, a convex cone  in the space of  linear operators on $A(X)$.\footnote{``$\mathrm{P}$'' is for Poincar\'e,  ``$\mathrm{L}$'' is for Lefschetz, and ``$\mathrm{K}$'' is for K\"ahler.}
For example,
 $A(X)$ may be  
 the cohomology of real $(q,q)$-forms on a compact K\"ahler manifold \cite{Gromov90}, 
the ring of algebraic cycles modulo homological equivalence on a smooth projective variety \cite{Grothendieck69},  
McMullen's algebra generated by the Minkowski summands of a simple convex polytope   \cite{McMullen93}, the combinatorial intersection cohomology of a convex polytope \cite{Karu04},
 the reduced Soergel bimodule of a Coxeter group element \cite{EW14}, or 
the Chow ring of a matroid defined in Section \ref{MatroidHR}.

Often, but not always, $A(X)$ has the structure of a graded algebra,  $\mathrm{P}$ is determined by the multiplicative structure of $A(X)$ up to a constant multiple, and $\mathrm{L}$ is the multiplication by an element in $A^1(X)$.
In any case, we expect the following properties to hold for the triple $(A(X),\mathrm{P}(X),\mathrm{K}(X))$ for every nonnegative integer $q \le \frac{d}{2}$:
\begin{enumerate}[(1)]
\item[(PD)] The bilinear pairing
\[
A^q(X) \times A^{d-q}(X) \longrightarrow \mathbb{R}, \qquad (\eta,\xi) \longmapsto \mathrm{P}(\eta,\xi)
\]
 is nondegenerate (the Poincar\'e duality for $X$).
\item[(HL)]  For any $\mathrm{L}_1,\ldots,\mathrm{L}_{d-2q} \in \mathrm{K}(X)$, the linear map
\[
A^q(X) \longrightarrow A^{d-q}(X), \qquad \eta \longmapsto \big(\prod_{i=1}^{d-2q} \mathrm{L}_i\big)\hspace{0.5mm} \eta
\] is bijective (the hard Lefschetz theorem for $X$).
\item[(HR)]  For  any $\mathrm{L}_0,\mathrm{L}_1,\ldots,\mathrm{L}_{d-2q} \in \mathrm{K}(X)$,  the bilinear form 
\[
A^q(X) \times A^q(X) \longrightarrow \mathbb{R}, \qquad (\eta_1,\eta_2) \longmapsto (-1)^q\ \mathrm{P}(\eta_1,\big(\prod_{i=1}^{d-2q} \mathrm{L}_i\big)\hspace{0.5mm} \eta_2)
\]
is positive definite on the kernel of the linear map
\[
A^q(X) \longrightarrow A^{d-q+1}(X), \qquad \eta \longmapsto \big(\prod_{i=0}^{d-2q} \mathrm{L}_i\big)\hspace{0.5mm} \eta
\]  
(the Hodge-Riemann relation for $X$).
\end{enumerate}
All three properties are known to hold for the objects listed above except one, which is the subject of  Grothendieck's standard conjectures on algebraic cycles.
The known proofs of the hard Lefschetz theorems and the Hodge-Riemann relations for different  objects 
 have certain structural similarities, but there is no known way of deducing one of them  from the others.


Hard Lefschetz theorems for various $X$'s have found numerous applications to problems of combinatorial nature. An early survey of these applications can be found in the proceedings article \cite{Stanley84}.
We highlight the following three:
\begin{enumerate}[(1)]
\item Erd\H{o}s-Moser conjecture \cite{Erdos65}, proved in \cite{Stanley80a}: Let $f(E,k)$ be the number of subsets of $E \subseteq \mathbb{R}$ whose elements sum to $k$. If
the cardinality of $E$ is $2n+1$, then
\[
f(E,k) \le f\Big([-n,n] \cap \mathbb{Z},0\Big).
\]
\item McMullen's $g$-conjecture  \cite{McMullen71}, proved in \cite{BL80,Stanley80b}:
The $f$-vector of a $d$-dimensional convex polytope $\mathrm{P}$ is the sequence $f_0(\mathrm{P}),\ldots,f_{d}(\mathrm{P})$,
where 
\[
f_i(\mathrm{P}) =\text{the number of $(i-1)$-dimensional faces of $\mathrm{P}$}.
\]
The $h$-vector of $\mathrm{P}$ is the sequence $h_0(\mathrm{P}),\ldots,h_{d}(\mathrm{P})$ defined by the identity
\[
 \sum_{i=0}^d h_i(\mathrm{P}) x^{i}=\sum_{i=0}^d f_i(\mathrm{P}) x^i (1-x)^{d-i}.
\]
The $g$-conjecture gives a complete numerical characterization of the $h$-vectors of simplicial polytopes.
In particular, for any $d$-dimensional simplicial polytope $\mathrm{P}$,
\[
h_i(\mathrm{P})=h_{d-i}(\mathrm{P}) \ \ \text{and} \ \ h_{i}(\mathrm{P}) \le h_{i+1}(\mathrm{P}) \ \ \text{for all} \ \  i<d/2.
\]
\item Dowling-Wilson conjecture \cite{DW74,DW75}, proved in \cite{HW}: 
Let $E$ be a finite subset of a vector space, and let $w_i(E)$ be the number of $i$-dimensional subspaces spanned by subsets of $E$. If $E$ spans a $d$-dimensional subspace, then
\[
w_i(E)\le w_{d-i}(E) \ \ \text{and} \ \ w_{i}(E) \le w_{i+1}(E) \ \ \text{for all} \ \  i<d/2.
\]
\end{enumerate}
None of the above statements has a known proof not using some version of $\text{HL}$.

When the Poincar\'e duality for $X$ is known, the Hodge-Riemann relation for $X$ is stronger than the hard Lefschetz theorem for $X$ in the sense that, for every $q$, there is a formal implication
\[
\text{$\text{HR}$ in degrees at most $q$} \Longrightarrow \text{$\text{HL}$ in degrees at most $q$}.
\]
In the remainder of this survey, we give an overview of applications of the Hodge-Riemann relations to concrete problems.
We remark that most known applications only use the following immediate consequence of $\text{HR}$ in degrees $q \le 1$: 
For any $\mathrm{L}_1,\ldots,\mathrm{L}_{d-2} \in \mathrm{K}(X)$, any matrix representing  the symmetric bilinear form
\[
A^1(X) \times A^1(X) \longrightarrow \mathbb{R}, \qquad (\eta_1,\eta_2) \longmapsto \mathrm{P}(\eta_1,\big(\prod_{i=0}^{d-2} \mathrm{L}_i\big)\eta_2)
\]
has exactly one positive eigenvalue.
One notable exception is the   implication
\begin{multline*}
\text{Grothendieck standard conjectures on algebraic cycles}\Longrightarrow \\
\text{Weil conjectures on zeta functions over finite fields},
\end{multline*}
which was one of the main motivations for formulating the standard conjectures  \cite{Corr01,Kle68,Kle94}.
It will be interesting to find applications of $\text{HR}$ for $q>1$ in other contexts too.









\section{Applications of the Hodge-Riemann relations}

\subsection{Mixed discriminants and permanents}\label{MD}

The notion of mixed discriminant arises when one combines the determinant with the matrix sum.
To define the mixed discriminant, let $\mathbf{A}=(\mathrm{A}_1,\ldots,\mathrm{A}_d)$ be a collection of real symmetric $d \times d$ matrices, 
and consider the  function
\[
\text{det}_\mathbf{A}: \mathbb{R}^d \longrightarrow \mathbb{R}, \qquad (t_1,\ldots,t_d)\longmapsto \text{det}(t_1 \mathrm{A}_1+\cdots+t_d \mathrm{A}_d),
\]
which is a homogeneous polynomial of degree $d$. The number
\[
D(\mathrm{A}_1,\ldots,\mathrm{A}_d)=\frac{\partial^d}{\partial t_1 \cdots \partial t_d}\text{det}_\mathbf{A}(0,\ldots,0)
\]
is called the \emph{mixed discriminant} of $\mathbf{A}$. The mixed discriminant is symmetric in $\mathbf{A}$, and it is nonnegative whenever all the matrices in $\mathbf{A}$ are positive semidefinite.\footnote{The latter fact can be viewed as a Hodge-Riemann relation in degree $0$.}

Now let $\mathbf{P}=(\mathrm{P}_1,\ldots,\mathrm{P}_{d-2})$ be any collection of  $d \times d$ positive semidefinite matrices.
Define a  symmetric bilinear form $\text{HR}(\mathbf{P})$ on the space of real symmetric $d \times d$ matrices by
\[
\text{HR}(\mathbf{P}):\text{Sym}_d \times \text{Sym}_d \longrightarrow \mathbb{R}, \qquad (\eta_1,\eta_2) \longmapsto D(\eta_1,\eta_2,\mathrm{P}_1,\ldots,\mathrm{P}_{d-2}).
\]
In \cite{Aleksandrov38}, Aleksandrov proved the following statement and used it in his proof of the Aleksandrov-Fenchel inequality for mixed volumes of convex bodies.
To avoid trivialities, we suppose that $\text{HR}(\mathbf{P})$ is not identically zero.

\begin{theorem}\label{HRDiscriminant}
Any matrix representing $\mathrm{HR}(\mathbf{P})$ has exactly one positive eigenvalue.
\end{theorem}

It follows from Cauchy's eigenvalue interlacing theorem that, for any positive semidefinite $d \times d$ matrices $\mathrm{A}_1,\ldots,\mathrm{A}_d$, 
\[
\text{det}\left(\begin{array}{cc}
D(\mathrm{A}_1,\mathrm{A}_1,\mathrm{A}_3,\ldots,\mathrm{A}_d) & D(\mathrm{A}_1,\mathrm{A}_2,\mathrm{A}_3,\ldots,\mathrm{A}_d) \\
D(\mathrm{A}_1,\mathrm{A}_2,\mathrm{A}_3,\ldots,\mathrm{A}_d) & D(\mathrm{A}_2,\mathrm{A}_2,\mathrm{A}_3,\ldots,\mathrm{A}_d)
\end{array}\right) \le 0.
\]
Theorem \ref{HRDiscriminant} is, in fact, a Hodge-Riemann relation in degree $1$.
The object $X$ is the $d$-dimensional complex vector space $\mathbb{C}^d$,
the algebra $A(X)$ is the ring of real differential forms with constant coefficients on $\mathbb{C}^d$,
and the cone $\mathrm{K}(X)$ is the spectrahedral cone of all $d \times d$ positive definite matrices.
Elementary proofs of the Hodge-Riemann relation for $X$ in any degree can be found in \cite{Gromov90,Timorin98}.

In the important special case  when all the matrices are diagonal, the mixed discriminant is a permanent.
Precisely,
if $\mathrm{A}=(a_{ij})$ is an $d \times d$ matrix and if $\mathrm{A}_i$ is the diagonal matrix whose $j$-th diagonal element is $a_{ij}$, then
\[
d! \hspace{0.5mm} D(\mathrm{A}_1,\ldots,\mathrm{A}_d)=\text{per}(\mathrm{A})\vcentcolon=\sum_{\sigma} \prod_{i=1}^d a_{i \sigma(i)},
\]
where $\sigma$ runs through all permutations of $\{1,\ldots,d\}$. 
Therefore, for any column vectors $a_1,\ldots,a_d$ in $\mathbb{R}^n$ with nonnegative entries, 
\[
\text{per}(a_1,a_2,a_3,\ldots,a_d)^2 \ge \text{per}(a_1,a_1,a_3,\ldots,a_d) \text{per}(a_2,a_2,a_3,\ldots,a_d).
\]
The above special case of the Hodge-Riemann relations for $\mathbb{C}^d$ was the main ingredient in Egorychev's and Falikman's proofs of van der Waerden's conjecture that
the permanent of any  doubly stochastic $d \times d$ matrix is at least $d!/d^d$.
See \cite{Knuth81,vanLint82} for more on van der Waerden's permanent conjecture.

\subsection{Mixed volumes of convex bodies}\label{MV}

The notion of mixed volume arises when one combines the volume with the Minkowski sum.
For any collection of convex bodies $\mathbf{P}=(\mathrm{P}_1,\ldots,\mathrm{P}_d)$ in $\mathbb{R}^d$, consider the  function
\[
\text{vol}_\mathbf{P}: \mathbb{R}^d_{\ge 0} \longrightarrow \mathbb{R}_{\ge 0},\qquad (t_1,\ldots,t_d)\longmapsto \text{vol}(t_1 \mathrm{P}_1+\cdots+t_d \mathrm{P}_d).
\]
Minkowski noticed that $\text{vol}_\mathbf{P}$ is  a homogeneous polynomial of degree $d$, and called the number
\[
V(\mathrm{P}_1,\ldots,\mathrm{P}_d)=\frac{\partial^d}{\partial t_1 \cdots \partial t_d}\text{vol}_\mathbf{P}(0,\ldots,0)
\]
 the \emph{mixed volume} of $\mathbf{P}$.
 The mixed volume is symmetric in $\mathbf{P}$, and it is nonnegative for any $\mathrm{P}$.\footnote{The latter fact can be viewed as a Hodge-Riemann relation in degree $0$.}

Now let $\eta_1,\ldots,\eta_n$ be another collection of convex bodies in $\mathbb{R}^d$, and define an $n \times n$ matrix $\textrm{AF}=(\textrm{AF}_{ij})$ by 
\[
\textrm{AF}_{ij}=V(\eta_i,\eta_j,\mathrm{P}_1,\ldots,\mathrm{P}_{d-2}).
\]
If   $\textrm{AF}\neq 0$, then the mixed volume analog of Theorem \ref{HRDiscriminant} holds.

\begin{theorem}\label{ConvexHR}
The matrix $\mathrm{AF}$ has exactly one positive eigenvalue.
\end{theorem}

It follows that the mixed volume satisfies the \emph{Aleksandrov-Fenchel inequality}
\[
\text{det}\left(\begin{array}{cc}
V(\mathrm{P}_1,\mathrm{P}_1,\mathrm{P}_3,\ldots,\mathrm{P}_d)& V(\mathrm{P}_1,\mathrm{P}_2,\mathrm{P}_3,\ldots,\mathrm{P}_d)\\
V(\mathrm{P}_1,\mathrm{P}_2,\mathrm{P}_3,\ldots,\mathrm{P}_d) & V(\mathrm{P}_2,\mathrm{P}_2,\mathrm{P}_3,\ldots,\mathrm{P}_d)
\end{array}\right) \le 0.
\]
In particular,
the sequence of mixed volumes of two convex bodies is log-concave:
\[
V(\underbrace{\mathrm{P}_1,\ldots,\mathrm{P}_1}_{i},\underbrace{\mathrm{P}_2,\ldots,\mathrm{P}_2}_{d-i})^2 \ge V(\underbrace{\mathrm{P}_1,\ldots,\mathrm{P}_1}_{i-1},\underbrace{\mathrm{P}_2,\ldots,\mathrm{P}_2}_{d-i+1})
V(\underbrace{\mathrm{P}_1,\ldots,\mathrm{P}_1}_{i+1},\underbrace{\mathrm{P}_2,\ldots,\mathrm{P}_2}_{d-i-1}).
\]
Aleksandrov reduced Theorem \ref{ConvexHR}  to the case when the Minkowski sum of all the relevant convex bodies, say $\mathrm{P}$, is a simple convex polytope.
Under this hypothesis, Theorem \ref{ConvexHR} is a Hodge-Riemann relation in degree $1$ \cite{Gromov90,McMullen93,Teissier79}.
The object $X$  is the convex polytope $P$,
the algebra $A(X)$ is McMullen's polytope algebra generated by the Minkowski summands of $P$, and the cone $\mathrm{K}(X)$ is the cone of convex polytopes that share the normal fan with $\mathrm{P}$.
Elementary proofs of the Hodge-Riemann relation for $X$ in any degree can be found in \cite{FK10,McMullen93,Timorin99}.






The Alexandrov-Fenchel inequality has been used to understand linear extensions of  partially ordered sets.
For example,  Chung, Fishburn, and Graham conjectured in \cite{CFG80} that, for any finite poset $\mathrm{Q}$,
\[
\text{Pr}_i(x)^2 \ge \text{Pr}_{i-1}(x)\text{Pr}_{i+1}(x) \ \ \text{for all $i$ and all $x \in \mathrm{Q}$},
\]
where $\text{Pr}_i(x)$ is the fraction of linear extensions of $\mathrm{Q}$ in which $x$ is the $i$-th largest element.
Stanley proved the conjecture by constructing suitable convex polytopes from  $x \in \mathrm{Q}$ and using the Alexandrov-Fenchel inequality \cite{Stanley81}.
Now let  $\text{Pr}(x_1<x_2)$ be the fraction of linear extensions of $\mathrm{Q}$ in which $x_1$ is smaller than $x_2$.
In \cite{KS84}, Kahn and Saks employed Stanley's method to deduce the following remarkable fact from the Alexandrov-Fenchel inequality:
\[
\text{If $\mathrm{Q}$ is not a chain, then there are elements $x_1,x_2 \in \mathrm{Q}$ such that $3/11< \text{Pr}(x_1<x_2) < 8/11$.}
\]  
This confirmed a conjecture of Fredman \cite{Fredman76} and Linial \cite{Linial84} that the information theoretic lower bound for the general sorting problem  is tight up to a multiplicative constant.

\subsection{The correlation constant of a field}

Let $G$ be a finite connected graph, let $i$, $j$ be distinct edges, and let $T$ be a random spanning tree of $G$.
Kirchhoff's effective resistance formula can be used to show that the probability that $i$ is in $T$ can only decrease by assuming that $j$ is in $T$: 
\[
\text{Pr}(i \in T) \ge \text{Pr}(i \in T \mid j \in T).
\]
In other words, the number $b_{-}$ of spanning trees containing given edges satisfies
\[
\frac{b_i}{b} \ge \frac{b_{ij}}{b_j}.
\]

Now let $\mathrm{M}$ be a finite spanning subset of a vector space $V$, let $i$, $j$ be distinct nonzero vectors in $\mathrm{M}$, and write $b_{-}$ for the number of bases in $\mathrm{M}$ containing given vectors.
Do we still have
\[
\frac{b_i}{b} \ge \frac{b_{ij}}{b_j}?
\]
In \cite{SW75}, Seymour and Welsh gave the first example of $\mathrm{M}$ over a field of characteristic $2$ with
$
\frac{b_{\hspace{0.5mm}}b_{ij}}{b_ib_j} = \frac{36}{35}$ for some $i$ and $j$.
How large can the ratio be?

\begin{definition}
The \emph{correlation constant} of a field $k$ is the supremum of $\frac{b_{\hspace{0.5mm}}b_{ij}}{b_ib_j}$ over all pairs of distinct vectors $i$ and $j$ in finite vector configurations in vector spaces over $k$.
\end{definition}

The correlation constant may be an interesting invariant of a field, although it is not immediately clear that the constant is finite. 
In fact, the finiteness of the correlation constant is one of the consequences of the Hodge-Riemann relations for vector configurations.
Let $n$ be the number of vectors in $\mathrm{M}$,  and let $\textrm{HR}(\mathrm{M})$ be the symmetric $n \times n$ matrix  
\[
\arraycolsep=1.1pt\def\arraystretch{1.3}
\textrm{HR}(\mathrm{M})_{ij} = \left\{\begin{array}{cl} 0 & \quad \text{if $i=j$,}\\ b_{ij} &\quad \text{if $i \neq j$.}\end{array}\right.
\]
To avoid the trivial case $\textrm{HR}(\mathrm{M})=0$, we suppose that the dimension of $V$ is at least $2$.
For example, if $\mathrm{K}_4$ is the set of six column vectors of the matrix
\[
\tiny
\left( \begin{array}{rrrrrr}
 1& 1& 1& 0& 0 &0\\
 -1& 0& 0& 1& 1 &0\\
 0& -1& 0& -1& 0 &1\\
 0& 0& -1& 0& -1&-1
 \end{array}\right),
\]
then $\textrm{HR}(\textrm{K}_4)$ is the $6 \times 6$ symmetric matrix
\[
\tiny
\left(\begin{array}{cccccc}
0 & 3 & 3 & 3 & 3 & 4\\
3 & 0 & 3 &3 & 4 & 3\\
3 & 3 & 0 &4 & 3 & 3\\
3 & 3 & 4 &0 & 3 & 3\\
3 & 4 & 3 &3 & 0 & 3\\
4 & 3 & 3 &3 & 3 & 0
\end{array}\right).
\]
In \cite{HW}, the following statement was deduced from Theorem \ref{MainTheorem}.

\begin{theorem}
The matrix $\textrm{HR}(\mathrm{M})$ has exactly one positive eigenvalue.
\end{theorem}

In fact, the same statement holds more generally for any matroid $\mathrm{M}$ \cite[Remark 15]{HW}. To deduce a bound on the correlation constant, consider the restriction of $\textrm{HR}(\mathrm{M})$ to the three-dimensional subspace of $\mathbb{R}^n$ spanned by $\mathbf{e}_i$, $\mathbf{e}_j$, and $(1,\ldots,1)$. 
Cauchy's eigenvalue interlacing theorem shows that the resulting $3 \times 3$ symmetric matrix also has exactly one positive eigenvalue.
Expressing the  $3 \times 3$  determinant, which should be nonnegative, we get the inequality
\[
\frac{b_{\hspace{0.5mm}}b_{ij}}{b_ib_j} \le  2-2(\text{dim} V)^{-1}.
\]
Thus the correlation constant of any field is at most $2$.
What is the correlation constant of, say,  $\mathbb{Z}/2\mathbb{Z}$? 
Does the correlation constant  depend on the field?

\subsection{The chromatic polynomial of a graph}\label{Chromatic}

Generalizing earlier work of Birkhoff,   Whitney introduced in \cite{Whitney32} the \emph{chromatic polynomial} of a connected graph $G$ as the function on $\mathbb{N}$ defined by
\[
\chi_G(q) = \text{the number of proper $q$-colorings of $G$}.
\]
In other words, $\chi_G(q)$ is the number of ways to color the vertices of $G$ using $q$ colors so that the endpoints of  every edge have different colors.
Whitney noticed that the chromatic polynomial is indeed a polynomial. 
In fact, we can write
\[
\chi_G(q)/q=a_0(G)q^{d}-a_1(G)q^{d-1}+\dots+(-1)^da_d(G)
\]
for some positive integers $a_0(G),\ldots,a_d(G)$,
where $d$ is one less than the number of vertices.

\begin{example}\label{SquareGraph}
The cycle $C_4$ with $4$ vertices and $4$ edges
has the chromatic polynomial
\[
\chi_{C_4}(q)={1}q^4-{4}q^3+{6}q^2-{3}q.
\]
\end{example}

The chromatic polynomial was originally devised as a tool for attacking the Four Color Problem, but soon it  attracted attention in its own right. 
In \cite{Read68}, Read conjectured  that
the coefficients of the chromatic polynomial form a unimodal sequence for any  graph. A few years later, Hoggar conjectured more generally in \cite{Hoggar74} that the coefficients form a log-concave sequence:
\[
a_i(G)^2 \ge a_{i-1}(G) a_{i+1}(G) \ \ \text{for any $i$ and $G$.}
\]
Notice that the chromatic polynomial can be computed using the \emph{deletion-contraction relation}:
if $G \backslash e$ is the deletion of an edge $e$ from $G$ and $G/e$ is the contraction of the same edge, then
\[
\chi_G(q)\ =\ \chi_{G\backslash e}(q) - \chi_{G/ e}(q).
\]
The first term counts the proper colorings of $G$, the second term counts the otherwise-proper colorings of $G$ where the endpoints of $e$ are permitted to have the same color, and the third term  counts the otherwise-proper colorings of $G$ where the endpoints of $e$ are mandated to have the same color.
For example, to compute the chromatic polynomial of the cycle $C_4$ in Example \ref{SquareGraph},  we write
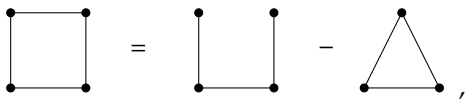
\begin{center}
\begin{tikzpicture}
\draw (0,0) -- (1,0);
\draw (0,0) -- (0,1);
\draw (1,1) -- (1,0);
\draw (1,1) -- (0,1);
\draw [fill=black] (0,0) circle (1.5pt); 
\draw [fill=black] (1,0) circle (1.5pt); 
\draw [fill=black] (0,1) circle (1.5pt); 
\draw [fill=black] (1,1) circle (1.5pt); 
\node at (1.7,0.5) {=};
\draw (2.5,0) -- (3.5,0);
\draw (2.5,0) -- (2.5,1);
\draw (3.5,1) -- (3.5,0);
\draw [fill=black] (2.5,0) circle (1.5pt); 
\draw [fill=black] (3.5,0) circle (1.5pt); 
\draw [fill=black] (2.5,1) circle (1.5pt); 
\draw [fill=black] (3.5,1) circle (1.5pt); 
\node at (4.2,0.5) {--};
\draw (4.7,0) -- (5.7,0);
\draw (5.2,1) -- (4.7,0);
\draw (5.7,0) -- (5.2,1);
\draw [fill=black] (4.7,0) circle (1.5pt); 
\draw [fill=black] (5.7,0) circle (1.5pt); 
\draw [fill=black] (5.2,1) circle (1.5pt); 
\node at (6.0,-0.1) {,};
\end{tikzpicture}
\end{center}
and use that the chromatic polynomials of the two smaller graphs are $q(q-1)^3$ and $q(q-1)(q-2)$, respectively. 
Note that, in general, the sum of  log-concave sequences need not be log-concave. 

The log-concavity conjecture for chromatic polynomials was proved in \cite{Huh12} by showing that the absolute values of the coefficients of $\chi_G(q)/(q-1)$ are mixed multiplicities of certain homogeneous ideals constructed from $G$.
The notion of mixed multiplicities is a commutative algebraic analog of the notion of  mixed volumes, and it can be shown that mixed multiplicities of homogeneous ideals satisfy a version of the Aleksandrov-Fenchel inequality.
To formulate the underlying Hodge-Riemann relation  in purely combinatorial terms was the primary motivation for \cite{AHK}.
The main result of \cite{AHK} will be reviewed in Section \ref{MatroidHR} below.

\subsection{Counting independent sets}\label{Independent}

How many linearly independent collection of $i$ vectors are there in a given configuration of vectors?
Let's write $\mathrm{M}$ for a finite subset of a vector space and $f_i(\mathrm{M})$ for the 
 number of independent subsets of $\mathrm{M}$ of size $i$.

\begin{example}\label{Fano}
Let $\mathrm{F}$ be the set of all nonzero vectors in the three-dimensional vector space over the field with two elements.
Nontrivial dependencies between elements of $\mathrm{F}$ can be read off from the picture of the Fano plane shown below.
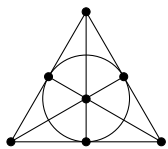
\begin{center}
\begin{tikzpicture}
\draw (0,0) -- (2,0);
\draw (2,0) -- (1,1.73);
\draw (1,1.73) -- (0,0);
\draw (0,0) -- (1.5,.866);
\draw (2,0) -- (.5,.866);
\draw (1,1.73) -- (1,0);
\draw (1,0.577) circle [radius=0.577];
\draw [fill=black] (0,0) circle (1.5pt); 
\draw [fill=black] (2,0) circle (1.5pt); 
\draw [fill=black] (1,1.73) circle (1.5pt); 
\draw [fill=black] (1.5,.866) circle (1.5pt); 
\draw [fill=black] (.5,.866) circle (1.5pt); 
\draw [fill=black] (1,0) circle (1.5pt); 
\draw [fill=black] (1,0.57) circle (1.5pt); 
\end{tikzpicture}
\end{center}
The nonempty independent subsets of $\mathrm{F}$ correspond to the seven points in $\mathrm{F}$, the twenty-one pairs of points in $\mathrm{F}$, and the twenty-eight triple of points in $\mathrm{F}$ not in one of the seven lines:
\[
{f_0(\mathrm{F})={1}, \quad f_1(\mathrm{F})={7}, \quad f_2(\mathrm{F})={21}, \quad f_3(\mathrm{F})={28}.}
\]
\end{example}

Welsh conjectured in 1971 that
the sequence $f_i(\mathrm{M})$ is unimodal  for any  $\mathrm{M}$ \cite{Welsh71}. Shortly after that, Mason conjectured more generally in \cite{Mason72} that the sequence is log-concave:
\[
f_i(\mathrm{M})^2 \ge f_{i-1}(\mathrm{M}) f_{i+1}(\mathrm{M}) \ \ \text{for any $i$ and $\mathrm{M}$.}
\]
In any small specific case, the conjecture can be verified by computing the $f_i(\mathrm{M})$'s by the \emph{deletion-contraction relation}:
if $\mathrm{M} \backslash v$ is the deletion of a nonzero vector $v$ from  $\mathrm{M}$ and $\mathrm{M}/v$ is the projection of $\mathrm{M}$ in the direction of $v$, then
\[
f_i(\mathrm{M})\ =\ f_{i}(\mathrm{M}\backslash v)+f_{i-1}(\mathrm{M}/v).
\]
The first term counts the number of independent subsets of size $i$, the second term counts the independent subsets of size $i$ not containing $v$,
and the third term counts the independent subsets of size $i$ containing $v$.
As in the case of graphs, we notice the apparent interference between the log-concavity conjecture and the additive nature of $f_i(\mathrm{M})$.

The log-concavity conjecture for $f_i(\mathrm{M})$ was proved in \cite{Lenz13}  by combining a geometric construction of \cite{HK12} and a matroid-theoretic construction of Brylawski \cite{Brylawski77}.
Given a spanning subset $\mathrm{M}$ of a $d$-dimensional vector space over a field $k$, one can construct a $d$-dimensional smooth projective variety $X(\mathrm{M})$ over $k$ and globally generated line bundles $L_1,L_2$ on $X(\mathrm{M})$ so that
\[
f_i(\mathrm{M})=\int_{X(\mathrm{M})}L_1^{d-i} L_2^i.
\]
The Hodge-Riemann relation for smooth projective varieties is known to hold in degrees $q \le 1$ \cite{Grothendieck58,Segre37}, and this implies  the log-concavity of $f_i(\mathrm{M})$ as in Sections \ref{MD}, \ref{MV}.
To express and verify the general Hodge-Riemann relation for $X(\mathrm{M})$  in purely combinatorial terms was another motivation for \cite{AHK}.

\subsection{The Hodge-Riemann relations for matroids}\label{MatroidHR}

In the 1930s, Hassler Whitney observed that several  notions in graph theory and linear algebra 
 fit together in a common framework, that of \emph{matroids} \cite{Whitney35}. This observation started  a new subject  with applications to a wide range of topics like characteristic classes, optimization, and moduli spaces. 
 
 \begin{definition}\label{FlatDef}
A \emph{matroid} $\mathrm{M}$ on a finite set $E$ is a collection of subsets of $E$, called {\emph{flats}} of $\mathrm{M}$, satisfying the following axioms: 
\begin{enumerate}[(1)]
\item The ground set $E$  is a flat.
\item If $F_1$ and $F_2$ are flats, then $F_1 \cap F_2$  is a flat.
\item If $F$ is a flat, then any element not in $F$ is contained in exactly one flat covering $F$.
\end{enumerate}
Here, a flat  is said to \emph{cover} another flat $F$ if it is minimal among the flats  properly containing $F$. 
\end{definition}

For our purposes, we may and will suppose that $\mathrm{M}$ is \emph{loopless}:
\begin{enumerate}[(1)]
\item[(4)] The empty subset of $E$  is a flat.
\end{enumerate}
Every maximal chain of flats in $F$ has the same length, and this common length is called the \emph{rank} of the flat $F$.
The rank of the flat $E$ is called the rank of the matroid $\mathrm{M}$.
Matroids are determined by their \emph{independent sets} (the idea of ``\emph{general position}''),
and can be alternatively defined in terms of independent sets \cite[Chapter 1]{Oxley11}.

\begin{example}
Let $E$ be the set of edges of a finite graph $G$. Call a subset $F$ of $E$ a flat when there is no edge in $E \setminus F$ whose endpoints are connected by a path in $F$. This defines a \emph{graphic matroid} on $E$.
\end{example}

\begin{example}
A \emph{projective space} $\mathbb{P}$ is a set with distinguished subsets, called \emph{lines}, satisfying:
\begin{enumerate}[(1)]
\item Every line contains more than two points.
\item If $x,y$ are distinct points, then there is exactly one line $xy$ containing $x$ and $y$.
\item If $x,y,z,w$ are distinct points, no three collinear, then 
\[ 
\text{the line $xy$ intersects the line $zw$} \Longrightarrow
\text{the line $xz$ intersects the line $yw$}.
 \]
\end{enumerate}
A subspace of $\mathbb{P}$ is a subset $\mathbb{S}$ of $\mathbb{P}$ such that 
\[
\text{$x$ and $y$ are distinct points in $\mathbb{S}$} \Longrightarrow \text{the line $xy$ is in $\mathbb{S}$}.
\]
For any finite subset $E$ of $\mathbb{P}$, the collection of sets of the form $E \cap \mathbb{S}$ has the structure of a matroid.
 Matroids arising from subsets of projective spaces over a field $k$ are said to be \emph{realizable} over $k$
(the idea of ``\emph{coordinates}'').
\end{example}


Not surprisingly, the notion of realizability is sensitive to the  field $k$. A matroid may arise from a vector configuration over one field while no such vector configuration exists over another field. 
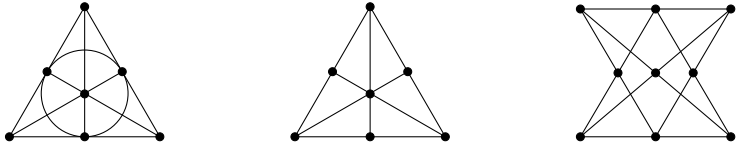
\begin{center}
\begin{tikzpicture}
\draw (0,0) -- (2,0);
\draw (2,0) -- (1,1.73);
\draw (1,1.73) -- (0,0);
\draw (0,0) -- (1.5,.866);
\draw (2,0) -- (.5,.866);
\draw (1,1.73) -- (1,0);
\draw (1,0.577) circle [radius=0.577];
\draw [fill=black] (0,0) circle (1.5pt); 
\draw [fill=black] (2,0) circle (1.5pt); 
\draw [fill=black] (1,1.73) circle (1.5pt); 
\draw [fill=black] (1.5,.866) circle (1.5pt); 
\draw [fill=black] (.5,.866) circle (1.5pt); 
\draw [fill=black] (1,0) circle (1.5pt); 
\draw [fill=black] (1,0.57) circle (1.5pt); 
\end{tikzpicture}
 \hspace{1.5cm}
\begin{tikzpicture}
\draw (0,0) -- (2,0);
\draw (2,0) -- (1,1.73);
\draw (1,1.73) -- (0,0);
\draw (0,0) -- (1.5,.866);
\draw (2,0) -- (.5,.866);
\draw (1,1.73) -- (1,0);
\draw [fill=black] (0,0) circle (1.5pt); 
\draw [fill=black] (2,0) circle (1.5pt); 
\draw [fill=black] (1,1.73) circle (1.5pt); 
\draw [fill=black] (1.5,.866) circle (1.5pt); 
\draw [fill=black] (.5,.866) circle (1.5pt); 
\draw [fill=black] (1,0) circle (1.5pt); 
\draw [fill=black] (1,0.57) circle (1.5pt); 
\end{tikzpicture}
\hspace{1.5cm}
\begin{tikzpicture}
\draw (0,0) -- (2,0);
\draw (0,1.7) -- (2,1.7);
\draw (0,0) -- (1,1.7);
\draw (0,0) -- (2,1.7);
\draw (1,0) -- (0,1.7);
\draw (1,0) -- (2,1.7);
\draw (2,0) -- (0,1.7);
\draw (2,0) -- (1,1.7);
\draw [fill=black] (0,0) circle (1.5pt); 
\draw [fill=black] (1,0) circle (1.5pt); 
\draw [fill=black] (2,0) circle (1.5pt); 
\draw [fill=black] (0.5,.85) circle (1.5pt); 
\draw [fill=black] (1,.85) circle (1.5pt); 
\draw [fill=black] (1.5,.85) circle (1.5pt); 
\draw [fill=black] (0,1.7) circle (1.5pt); 
\draw [fill=black] (1,1.7) circle (1.5pt); 
\draw [fill=black] (2,1.7) circle (1.5pt); 
\end{tikzpicture}
\end{center}
Among the rank $3$ matroids pictured above, where rank $1$ flats are represented by points and rank $2$ flats containing more than $2$ points are represented by lines, the first  is realizable over $k$ if and only if the characteristic of $k$ is $2$, the second  is realizable over $k$ if and only if the characteristic of $k$ is not $2$, and the third  is not realizable over any field.
It was recently shown that  almost all matroids are not realizable over any field \cite{Nelson}.

\begin{definition}\label{DefinitionMatroidChowRing}
We introduce variables $x_{F}$,
one for each nonempty proper flat $F$ of $\mathrm{M}$, and consider the polynomial ring
\[
S(\mathrm{M})=\mathbb{R}[x_{F}]_{F \neq \varnothing, F \neq E}.
\]
The \emph{Chow ring} $A(\mathrm{M})$ of $\mathrm{M}$ is the quotient of $S(\mathrm{M})$
by the ideal  generated by  the linear forms
\[
\sum_{i_1 \in F} x_{F} - \sum_{i_2 \in F} x_{F},
\]
one for each pair of distinct elements $i_1$ and $i_2$ of  $E$, and the quadratic monomials
\[
x_{F_1}x_{F_2},
\]
one for each pair of incomparable nonempty proper flats $F_1$ and $F_2$ of $\mathrm{M}$.
We have 
\[
A(\mathrm{M})=\bigoplus_{q} A^q(\mathrm{M}),
\]
where $A^q(\mathrm{M})$ is the span of degree $q$ monomials in $A(\mathrm{M})$.
\end{definition}

Feichtner and Yuzvinsky introduced the Chow ring of $\mathrm{M}$ in \cite{FY04}. When $\mathrm{M}$ is realizable over a field $k$, it is the Chow ring of the ``wonderful'' compactification of the complement of a hyperplane arrangement defined over $k$ studied by De Concini and Procesi \cite{DCP95}.

To formulate the hard Lefschetz theorem and the Hodge-Riemann relations for $A(\mathrm{M})$, we define a matroid analog of the K\"ahler cone in complex geometry.

\begin{definition}
A real-valued function $c$ on $2^E$  is said to be \emph{strictly submodular} if
\[
c_\varnothing=0, \quad c_E=0,
\]
and, for any two incomparable subsets $I_1,I_2 \subseteq E$, 
\[
 c_{I_1}+c_{I_2} > c_{I_1\,\cap\,I_2} +c_{I_1\,\cup\,I_2}.
 \]
 A strictly submodular function $c$ defines an element
\[
\mathrm{L}_c= \sum_F c_F x_F\in A^1(\mathrm{M}).
\]
The cone $\mathrm{K}(\mathrm{M})$ is defined to be the set of all such elements  in $A^1(\mathrm{M})$. 
\end{definition}

Now let $d+1$ be the rank of $\mathrm{M}$,
and write ``$\text{deg}$'' for the unique linear isomorphism
\[
\text{deg}\colon A^d(\mathrm{M}) \longrightarrow \mathbb{R} 
\]
which maps $x_{F_1}\cdots x_{F_d}$ to $1$ for every maximal chain 
$F_1 \subsetneq   \cdots \subsetneq F_d$ of nonempty proper flats \cite[Proposition 5.10]{AHK}.
We are ready to state the hard Lefschetz theorem and the Hodge-Riemann relation for $\mathrm{M}$ \cite[Theorem 8.9]{AHK}.

\begin{theorem}\label{MainTheorem}
Let $q$ be a nonnegative integer  $\le \frac{d}{2}$, and let $\mathrm{L}_0,\mathrm{L}_1,\ldots,\mathrm{L}_{d-2q} \in \mathrm{K}(\mathrm{M})$.
\begin{enumerate}[(1)]\itemsep 5pt
\item[(PD)] The product in $A(\mathrm{M})$ defines a nondegenerate bilinear pairing
\[
A^q(\mathrm{M}) \times A^{d-q}(\mathrm{M}) \longrightarrow \mathbb{R}, \qquad (\eta,\xi) \longmapsto \text{deg}(\eta \hspace{0.5mm}\xi).
\]
\item[(HL)] The multiplication by $\mathrm{L}_1,\ldots,\mathrm{L}_{d-2q}$ defines a linear bijection
\[
A^q(\mathrm{M}) \longrightarrow A^{d-q}(\mathrm{M}), \quad \eta \longmapsto  \big( \prod_{i=1}^{d-2q} \mathrm{L}_{i}\big) \hspace{0.5mm} \eta.
\]
\item[(HR)] The symmetric bilinear form 
\[
A^q(\mathrm{M}) \times A^{q}(\mathrm{M}) \longrightarrow \mathbb{R}, \qquad (\eta_1,\eta_2) \longmapsto (-1)^q\ \text{deg}\big( \big( \prod_{i=1}^{d-2q} \mathrm{L}_{i}\big)\hspace{0.5mm}\eta_1\eta_2\big)
\]
is positive definite on the kernel of the multiplication map
\[
A^q(\mathrm{M}) \longrightarrow A^{d-q+1}(\mathrm{M}), \quad \eta \longmapsto  \big( \prod_{i=0}^{d-2q} \mathrm{L}_{i}\big) \hspace{0.5mm} \eta.
\]
\end{enumerate}
\end{theorem}

We highlight the following consequence of $\textrm{HR}$ in degrees $\le 1$: 
For any $\xi_1, \xi_2 \in \mathrm{K}(\mathrm{M})$, 
\[
\left(\begin{array}{cc}
\text{deg}\big( \big( \prod_{i=1}^{d-2} \mathrm{L}_{i}\Big)\hspace{0.5mm}\xi_1\xi_1\big) & \text{deg}\big( \big( \prod_{i=1}^{d-2} \mathrm{L}_{i}\big)\hspace{0.5mm}\xi_1\xi_2\big) \\
\text{deg}\big( \big( \prod_{i=1}^{d-2} \mathrm{L}_{i}\Big)\hspace{0.5mm}\xi_1\xi_2\big) & \text{deg}\big( \big( \prod_{i=1}^{d-2} \mathrm{L}_{i}\big)\hspace{0.5mm}\xi_2\xi_2\big)
\end{array}\right)
\]
has exactly one positive eigenvalue.
Taking the determinant, we get an analog of the Alexandrov-Fenchel inequality
\[
 \text{deg}\big( \big( \prod_{i=1}^{d-2} \mathrm{L}_{i}\big)\hspace{0.5mm}\xi_1\xi_2\big)^2 \ge \text{deg}\big( \big( \prod_{i=1}^{d-2} \mathrm{L}_{i}\big)\hspace{0.5mm}\xi_1\xi_1\big)\text{deg}\big( \big( \prod_{i=1}^{d-2} \mathrm{L}_{i}\big)\hspace{0.5mm}\xi_2\xi_2\big).
\]
We apply the inequality to the \emph{characteristic polynomial} $\chi_{\mathrm{M}}(q)$,
 a generalization of the chromatic polynomial $\chi_G(q)$ to a matroid $\mathrm{M}$ that is not necessarily graphic \cite[Chapter 15]{Welsh76}. 
For this, we consider two distinguished elements of $A^1(\mathrm{M})$. For fixed $j \in E$, the elements are
\[
\alpha=\sum_{j \in F} x_F,\quad \beta=\sum_{j \notin F} x_F.
\]
The two elements do not depend on the choice of $j$, and they are limits of elements of the form $\ell_c$ for a strictly submodular function $c$. 
A bijective counting argument in \cite{AHK} shows that
\[
e_i(\mathrm{M})=\text{deg}(\alpha^i\ \beta^{d-i}) \ \ \text{for every $i$},
\]
where $e_i(\mathrm{M})$ is the sequence of integers satisfying the identity
\[
\chi_{\mathrm{M}}(q)/(q-1)
=e_0(\mathrm{M})q^d-e_1(\mathrm{M})q^{d-1}+\cdots+(-1)^d e_d(\mathrm{M}).
\]
Thus the sequence $e_i(\mathrm{M})$ is log-concave, which implies the following conjecture of Rota, Heron, and Welsh \cite{Heron72,Rota71,Welsh76}:
\[
\text{The  coefficients of  $\chi_\mathrm{M}(q)$ form a log-concave sequence for any matroid $\mathrm{M}$.}
\]
The above implies the  log-concavity of the sequence $a_i(G)$ in Section \ref{Chromatic} and the log-concavity of the sequence $f_i(\mathrm{M})$  in Section \ref{Independent}.
See  \cite[Chapter 15]{Oxley11} and \cite[Chapter 8]{White87}  for overviews and historical accounts.

\subsection{The reliability polynomial of a network}

Let $G$ be a finite connected graph with $v$ vertices and $n$ edges. 
The \emph{reliability} of $G$ is the probability that any two vertices remain connected when each edge is independently removed with the same probability $q$. 
Let's write $o_i(G)$ for the number of $i$-edge operational states.
For example, $o_{v-1}(G)$ is the number of spanning trees and $o_{n-1}(G)$ is the number of non-bridges.
Thus the reliability of $G$ is 
\[
R_G(q)=\sum_{i} o_i(G) (1-q)^{i} q^{n-i}.
\]
We define a sequence of integers $h_0(G),\ldots,h_d(G)$ by the identity
\[
R_G(q)/(1-q)^{v-1}= h_d(G) q^d+ h_{d-1}(G) q^{d-1}+\cdots+ h_0(G),
\]
where $d$ is one more than  the difference $n-v$.

\begin{example}
The complete graph on $4$ vertices has the reliability polynomial
\begin{align*}
R_{K_4}(q)=16 q^3(1-q)^3+ 15 q^2(1-q)^4+6q(1-q)^5 +1(1-q)^6
=(1-q)^3(6q^3+6q^2+3q+1).
\end{align*}
\end{example}

The numbers $h_i$ are closely related to the numbers $f_i$ of independent sets in Section \ref{Independent}.
Writing $\mathrm{M}$ for the dual of the graphic matroid of $G$, we have
\[
 \sum_{i=0}^d h_i(G) x^{i} =\sum_{i=0}^d f_i(\mathrm{M}) x^i (1-x)^{d-i}= \sum_{i=0}^d h_i(\mathrm{M}) x^{i}.
\]
Dawson conjectured in \cite{Dawson84} that the sequence $h_i(\mathrm{M})$ defined by the second equality is log-concave for any matroid $\mathrm{M}$:
\[
h_i(\mathrm{M})^2 \ge h_{i-1}(\mathrm{M})h_{i+1}(\mathrm{M}) \ \ \text{for any $i$ and $\mathrm{M}$.}
\]
Colbourn independently conjectured the same in the context of reliability polynomials \cite{Colbourn87}.

When $\mathrm{M}$ is the dual of a graphic matroid, or more generally when $\mathrm{M}$ is realizable over the complex numbers,
the log-concavity conjecture for $h_i(\mathrm{M})$ was proved in \cite{Huh15} by applying an algebraic analog of the Alexandrov-Fenchel inequality to the variety of critical points of the master function of a realization of $\mathrm{M}$ \cite{DGS12}. 
The underlying combinatorial Hodge-Riemann relation is yet to be formulated, and Dawson's conjecture for general matroids remains open.
The argument in the complex realizable case is tightly connected to the geometry of characteristic cycles \cite{Huh13}, suggesting that the combinatorial Hodge-Riemann relation  in this context will be strictly stronger than that of Section \ref{MatroidHR}.

\subsection{Unsolved problems}\label{Open}

The  log-concavity of a sequence is not only important because of its applications but because it hints the existence of a structure that satisfies $\text{PD}$, $\text{HL}$, and $\text{HR}$.
We close by listing some of the most interesting sequences that are conjectured to be log-concave.

\begin{enumerate}[(1)]
\item Rota's unimodality conjecture \cite{Rota71}:
If $w_k(\mathrm{M})$ is the number of rank $k$ flats of a rank $d$ matroid $\mathrm{M}$, then the sequence $w_0(\mathrm{M}),\ldots,w_d(\mathrm{M})$ is unimodal.
Welsh conjectured more generally in \cite{Welsh76} that the sequence is log-concave.
\item Fox's trapezoidal conjecture \cite{Fox62}: 
The sequence of absolute values of the coefficients of the Alexander polynomial of an alternating knot strictly increases, possibly plateaus, then strictly decreases.
Stoimenow conjectured more generally in \cite{Stoimenow05} that the sequence is log-concave.
\item Kazhdan-Lusztig polynomials of matroids \cite{EPW16}: For any matroid $\mathrm{M}$, the coefficients of the Kazhdan-Lusztig polynomial of $\mathrm{M}$ form a nonnegative log-concave sequence.
\end{enumerate}

\end{document}